\documentclass[12pt]{article}

\usepackage{amsmath, amsfonts, amssymb}
\usepackage{latexsym}
\usepackage{graphicx}
\usepackage{color}
\usepackage{url}

\newtheorem{thm}{Theorem}[section]
 \newtheorem{cor}[thm]{Corollary}

 \newtheorem{exa}[thm]{Example}

\newcommand{\iy}{\infty}

\newcommand{\bT}{{\mathbf T}}
\newcommand{\bZ}{{\mathbf Z}}

\newcommand{\cA}{{\mathcal A}}

\newcommand{\cL}{{\mathcal L}}

\newcommand{\xe}{\xi_\de\eta_\ga}

\newcommand{\al}{\alpha}

\newcommand{\ga}{\gamma}
\newcommand{\Ga}{\Gamma}
\newcommand{\de}{\delta}

\newcommand{\ka}{\kappa}
\newcommand{\la}{\lambda}

\newcommand{\ph}{\varphi}

\newcommand{\tht}{\theta}

\newcommand{\n}{\|}
\newcommand{\vsk}{\vspace{1mm}}
\newcommand{\vsg}{\vspace{2mm}}
\newcommand{\ti}{\widetilde}
\newcommand{\ov}{\overline}

\newcommand{\Ret}{{\rm Re}\,}

\renewcommand{\thefootnote}{\fnsymbol{footnote}}

\begin{document}

\begin{center}
{\Large \bf Toeplitz determinants with perturbations\\[0.5ex] in the corners}

\vspace{7mm}
{\large Albrecht B\"ottcher, Lenny Fukshansky,\\[0.5ex]

Stephan Ramon Garcia, Hiren Maharaj}
\end{center}

\medskip
\begin{quote}
\footnotesize{The paper is devoted to exact and asymptotic formulas for the determinants of Toeplitz matrices
with perturbations by blocks of fixed size in the four corners. If the norms of the inverses of the unperturbed
matrices remain bounded as the
matrix dimension goes to infinity, then standard perturbation theory yields asymptotic expressions for the
perturbed determinants. This premise is not satisfied for matrices generated by so-called Fisher-Hartwig symbols.
In that case we establish formulas for pure single Fisher-Hartwig singularities and for Hermitian matrices
induced by general Fisher-Hartwig symbols.
}
\let\thefootnote\relax\footnote{\hspace*{-7.5mm} MSC 2010: Primary 47B35; Secondary 15A15, 15B05}
\let\thefootnote\relax\footnote{\hspace*{-7.5mm} Keywords: Toeplitz matrix, Toeplitz determinant, Fisher-Hartwig symbol}
\let\thefootnote\relax\footnote{\hspace*{-7.5mm} Fukshansky acknowledges support by Simons Foundation grant \#279155,
Garcia acknowledges support by NSF grant DMS-1265973. }
\end{quote}

\section{Introduction}\label{S1}

This paper was prompted by a problem from lattices associated with finite Abelian groups.
This problem, which will be described in Section \ref{SL}, led to the computation
of the determinant of the $n \times n$ analogue $A_n$ of the matrix
\begin{equation}
A_6= \left(\begin{array}{rrrrrr}
6 & -4 & 1 & 0 & 0 & 1\\
-4 & 6 & -4 & 1 & 0 & 0\\
1 & -4 & 6 & -4 & 1 & 0\\
0 & 1 & -4 & 6 & -4 & 1\\
0 & 0 & 1 & -4 & 6 & -4\\
1 & 0 & 0 & 1 & -4 & 6
\end{array}\right).\label{A6}
\end{equation}
It turns out that $\det A_n =(n+1)^3$.
What makes the matter captivating is that the determinant of the $n \times n$ version $T_n$ of
\begin{equation}T_6= \left(\begin{array}{rrrrrr}
6 & -4 & 1 & 0 & 0 & 0\\
-4 & 6 & -4 & 1 & 0 & 0\\
1 & -4 & 6 & -4 & 1 & 0\\
0 & 1 & -4 & 6 & -4 & 1\\
0 & 0 & 1 & -4 & 6 & -4\\
0 & 0 & 0 & 1 & -4 & 6
\end{array}\right)\label{T6}
\end{equation}
is a so-called pure Fisher-Hartwig determinant. The latter determinant is known to be
\begin{equation}
\frac{(n+1)(n+2)^2(n+3)}{12}.\label{adet}
\end{equation}
This formula was established in \cite{BoSiJFA}. See also \cite[Theorem 10.59]{BoSiBu} or \cite{BoWi}.
We were intrigued by the question why the perturbations in the corners lower the growth from
$n^4$ to~$n^3$.

\vsk
The general context is as follows.
Every complex-valued function $a \in L^1$ on the unit circle $\bT$ has well-defined Fourier coefficients
\[a_k=\frac{1}{2\pi}\int_{-\pi}^\pi a(e^{i\tht}) e^{-ik\tht}\,d\tht, \quad k \in \bZ,\]
and generates the infinite Toeplitz matrix $T(a)=(a_{j-k})_{j,k=1}^\iy$. The principal $n \times n$
truncation of this matrix is denoted by $T_n(a)$. Thus, $T_n(a)=(a_{j-k})_{j,k=1}^n$.
The function $a$ is usually referred to as the symbol of the infinite matrix $T(a)$ and of the
sequence $\{T_n(a)\}_{n=1}^\iy$. For example, matrix (\ref{T6}) is just $T_6(a)$ with
\begin{equation}
a(t)=t^{-2}-4t^{-1}+6-4t+t^2=\left(1-\frac{1}{t}\right)^2(1-t)^2=|1-t|^4,\label{1.3}
\end{equation}
where here and in the following $t=e^{i\tht}$. The function $a(t)=|1-t|^4$ has a zero on the
unit circle and therefore the classical Szeg\H{o} limit theorem cannot be used to compute
$\det T_n(a)$ asymptotically. Fortunately, $a(t)=|1-t|^4$ is a special pure Fisher-Hartwig symbol,
and for such symbols the determinants are known both exactly and asymptotically.

\vsk
In Section \ref{S2} we consider the determinants of perturbations of $T_n(a)$ under the assumption
that the norms of the inverses of $T_n(a)$ remain bounded as $n \to \iy$. In that case, under
mild additional conditions, the determinants of the unperturbed matrices are asymptotically given
by Szeg\H{o}'s strong limit theorem.

\vsk
The (standard) techniques of Section~\ref{S2} do not work
for so-called Fisher-Hartwig symbols. This class of symbols was introduced by Fisher and Hartwig
in~\cite{FiHa} in connection with several problems of statistical physics. Paper~\cite{Dei}
contains a very readable exposition of the entire story up to the recent developments. See also
the books~\cite{BoSiUni} and~\cite{BoSiBu}. A pure Fisher-Hartwig symbol is of the form
$a(t)=(1-t)^\ga(1-1/t)^\de$. In particular, symbol~(\ref{1.3}) is of this form with $\ga=\de=2$.
Determinants of perturbed Toeplitz matrices with pure Fisher-Hartwig symbols are studied in
Section~\ref{S3}. Among other things, we there give an explanation of the growth drop from
$n^4$ to $n^3$ when replacing~(\ref{T6}) by~(\ref{A6}).

\vsk
In Section \ref{S4} we consider the very general case of symbols $a \in L^1$ which are
nonnegative a.e. on the unit circle and whose logarithm $\log a$ is also in $L^1$. We
there show that the quotient of the perturbed and unperturbed determinants approaches
a limit as $n \to \iy$ and we determine this limit. The class of symbols treated
in Section~\ref{S4} includes the general positive Fisher-Hartwig
symbols $a(t)=|t_1-t|^{2\al_1}\cdots |t_r-t|^{2\al_r} b(t)$ where the $t_j$
are distinct points on $\bT$, the $\al_j$ are real numbers in $(-1/2,1/2)$,
and $b$ is a sufficiently smooth and strictly positive function on $\bT$.

\section{The lattice of a cyclic group} \label{SL}

The idea behind paper \cite{FM} is to associate a lattice with an elliptic curve
and then to connect arithmetic properties of the curve with geometric properties of the lattice.
The lattices obtained in this way are generated by finite Abelian (additively written) groups
$G=\{0,P_1, \ldots, P_n\}$ and are of the form
\begin{equation}
\{(x_1, \ldots, x_n, -x_1-
\cdots-x_n)\in \bZ^{n+1}: x_1P_1+\cdots+x_n P_n=0\}.\label{LG}
\end{equation}
One may think of these lattices as  full rank sublattices of the well-known family of root lattices
\[\cA_{n}:=\{(x_1, \ldots, x_{n}, -x_1-\cdots-x_{n})\in \bZ^{n+1}: x_1,\ldots, x_n \in \bZ\}.\]
A fundamental quantity of every lattice is its determinant (i.e., the volume of a fundamental domain). Papers~\cite{BFGM} and~\cite{Sha}
contain a simple, purely group-theoretic
argument which shows that the determinant of the lattice~(\ref{LG}) equals $(n+1)^{3/2}$.
In particular, the determinant depends only on the order of the group. As shown in~\cite{BFGM},
this result can also be derived
in a completely elementary fashion via the computation of (usual) determinants. Here is this computation in the simple case
where $G$ is the cyclic group of order $n+1$. The corresponding lattice is
\[\cL_n:=\{(x_1, \ldots, x_n, -x_1\cdots-x_n)\in \bZ^{n+1}: x_1+2x_2+\cdots+n x_n=0\;\:\mbox{modulo}\;\:n+1\}.\]
The rank of the lattice $\cL_n \subset \cA_n$ is $n$, and in~\cite{BFGM} it is proved that the columns of the $(n+1) \times n$ matrix
\[B_n=\left(\begin{array}{rrrrrr}
-2 & 1 & 0 & \ldots  & 0 & 0\\
1 & -2 & 1 & \ldots  & 0 & 0\\
0 & 1 & -2 & \ddots & 0 & 0\\
\vdots & \ddots & \ddots & \ddots & \ddots & \vdots\\
0 & 0 & 0 & \ddots & -2 & 1\\
0 & 0 & 0 & \ldots & 1 & -2\\
1 & 0 & 0 & \ldots & 0 & 1
\end{array}\right)\]
form a basis of the lattice $\cL_n$. The determinant of $\cL_n$ is known to be $\sqrt{\det(B_n^\top B_n)}$,
and $B_n^\top B_n$ is just the matrix $A_n$ we encountered in the introduction.
Thus, the calculation of the determinant of the lattice $\cL_n$ is equivalent to
the computation of the determinant of the matrix $A_n$.

\vsk
Applying the Cauchy-Binet formula, we may write
\[\det A_n =\det(B_n^\top B_n)=(\det C_1)^2+(\det C_2)^2+\cdots+(\det C_{n+1})^2,\]
where $C_j$ results from $B_n$ by deleting the $j$th row.
Expanding $\det C_j$ along the last row and using the fact that
the determinant of the $k \times k$
tridiagonal Toeplitz matrix with $-2$ on the main diagonal and $1$ on the two neighboring diagonals
is $(-1)^k (k+1)$, it follows that each $\det C_j$ equals $\pm (n+1)$. Consequently,
\[\det A_n = (n+1)\cdot (n+1)^2=(n+1)^3,\]
as desired.

\section{The tame case}\label{S2}

We now turn to Toeplitz determinants and their perturbations.
Suppose the symbol $a$ is a piecewise continuous function, that is, the one-sided limits
$a(t\pm 0)$ exist for each $t \in \bT$. Let $a^\sharp$ be the continuous curve in the plane
that results from the range of $a$ by filling in the line segments $[a(t-0),a(t+0)]$ for each
$t$ where $a$ makes a jump. A famous theorem of Gohberg~\cite{Go} (see also~\cite[Corollary 2.19]{BoSiUni}
or~\cite[Theorem IV.4.1]{GoFe})
says that if the curve $a^\sharp$ does not pass through the origin and has winding number zero
about the origin, then the infinite matrix $T(a)$ generates
a bounded and invertible operator on $\ell^2$, the truncations $T_n(a)$ are invertible
for all sufficiently large $n$, and the inverses $T_n^{-1}(a):=[T_n(a)]^{-1}$ converge
strongly to the inverse $T^{-1}(a):=[T(a)]^{-1}$. To be more precise,
\begin{equation}
T_n^{-1}(a)P_n x \;\:\mbox{converges in}\;\: \ell^2\;\:\mbox{to}\;\:T^{-1}(a)x\;\:\mbox{for every}\;\: x \in \ell^2, \label{fsm}
\end{equation}
where $P_n$ is the projection $P_n :\{x_1, x_2, x_3, \ldots\} \mapsto \{x_1, \ldots, x_n, 0,\ldots\}$.

\vsk
Let $E_{11}, E_{12}, E_{21}, E_{22}$ be four $m_0 \times m_0$ matrices. For $n \ge 2m_0$, we denote
by $E_n$ the $n \times n$ matrix with the matrices $E_{jk}$ in the corners and zeros elsewhere,
\[E_n=\left(\begin{array}{ccc}
E_{11} & 0 & E_{12}\\
0 & 0 & 0\\
E_{21} & 0 & E_{22}\end{array}\right).\]
If $T(a)$ is invertible, then the operator $T^{-1}(a)$ is given by an infinite matrix
in the natural fashion.
We denote the entries of $T^{-1}(a)$ by $c_{jk}$ and let $S_{11}
=(c_{jk})_{j,k=1}^{m_0}$ stand for the upper-left $m_0 \times m_0$ block of $T^{-1}(a)$,
\[T^{-1}(a)=\left(\begin{array}{cccc}
c_{11} & \ldots & c_{1, m_0} & \ldots\\
\ldots & & \ldots & \ldots\\
c_{m_0, 1} & \ldots & c_{m_0,m_0} & \ldots\\
\ldots & \ldots & \ldots & \ldots \end{array}\right)
=\left(\begin{array}{cc} S_{11} & \ast\\
\ast & \ast\end{array}\right).\]
Let $W_m$ be the $m \times m$ counter-identity matrix,
that is, $W_m$ has ones on the counter-diagonal and zeros elsewhere. Given an $m \times m$ matrix $B$,
we denote by $\ti{B}$ the matrix $W_mBW_m$.  Recall that $B^\top$ stands for the transposed matrix. Toeplitz matrices
enjoy the property that $[T_n(a)]\,\ti{~}=[T_n(a)]^\top = T_n(\ti{a})$, where $\ti{a}$ is the function defined
by $\ti{a}(t)=a(1/t)$, $t \in \bT$.

\begin{thm} \label{Theo 2.1}
Let $a$ be piecewise continuous and suppose $a^\sharp$ does not contain the origin and has winding number
zero about the origin.
Then
\begin{equation}
\lim_{n\to \iy}\frac{\det (T_n(a)+E_n)}{\det T_n(a)}
=\det\left[\left(\begin{array}{cc}
I & 0\\ 0 & I\end{array}\right)+\left(\begin{array}{cc}
S_{11} & 0\\
0 & \ti{S}_{11}^\top
\end{array}\right)
\left(\begin{array}{cc}
E_{11} & E_{12}\\
E_{21} & E_{22}\end{array}\right)\right].\label{limit}
\end{equation}
\end{thm}

{\em Proof.} We know that $T_n(a)$ is invertible for sufficiently large $n$, in which case
\begin{equation}
\det (T_n(a)+E_n)=\det T_n(a) \det (I+T_n^{-1}(a)E_n). \label{2.1}
\end{equation}
We write $T_n^{-1}(a)$ as
\begin{equation}
T_n^{-1}(a)=\left(\begin{array}{ccc}
S_{11}^{(n)} & \ast & S_{12}^{(n)}\\
\ast & \ast & \ast\\
S_{21}^{(n)} & \ast & S_{22}^{(n)}\end{array}\right)\label{Sjk}
\end{equation}
with $m_0 \times m_0$ matrices $S_{jk}^{(n)}$. From (\ref{fsm}) we infer that if $I$ is the
$m_0 \times m_0$ identity matrix, then
\[\left(\begin{array}{c}
S_{11}^{(n)}\\
\ast\\
S_{21}^{(n)}\end{array}\right)
=T_n^{-1}(a)\left(\begin{array}{c}
I\\
0\\
0\end{array}\right) \to T^{-1}(a)\left(\begin{array}{c}
I \\ 0 \end{array}\right)
= \left(\begin{array}{c}
S_{11} \\ \ast \end{array}\right).\]
This implies that $S_{11}^{(n)} \to S_{11}$ and $S_{21}^{(n)} \to 0$.
(Here we are dealing with convergence of $m_0 \times m_0$ matrices, which may be understood
entry-wise.) We further have
\[T_n^{-1}(\ti{a})=W_n T_n^{-1}(a)W_n=\left(\begin{array}{ccc}
\ti{S}_{22}^{(n)} & \ast & \ti{S}_{21}^{(n)}\\
\ast & \ast & \ast\\
\ti{S}_{12}^{(n)} & \ast & \ti{S}_{11}^{(n)}\end{array}\right)\]
and
\[[T_n^{-1}(a)]^\top=\left(\begin{array}{ccc}
~[S_{11}^{(n)}]^\top & \ast & [S_{21}^{(n)}]^\top\\
~\ast & \ast & \ast\\
~[S_{12}^{(n)}]^\top & \ast & [S_{22}^{(n)}]^\top
\end{array}\right).
\]
Since $T_n^{-1}(\ti{a})=[T_n^{-1}(a)]^\top$, we see that $\ti{S}_{22}^{(n)}=[S_{11}^{(n)}]^\top$
and $\ti{S}_{12}^{(n)}=[S_{21}^{(n)}]^\top$. From what was already proved we therefore obtain that
$S_{12}^{(n)} \to 0$ and $S_{22}^{(n)}=[\ti{S}_{11}^{(n)}]^\top \to \ti{S}_{11}^\top.$
The matrix $I+T_n^{-1}(a)E_n$ equals
\[
\left(\begin{array}{ccc}
I+S_{11}^{(n)}E_{11}+S_{12}^{(n)}E_{21} & 0 & S_{11}^{(n)}E_{12}+S_{12}^{(n)}E_{22}\\
0 & I & 0\\
S_{21}^{(n)}E_{11}+S_{22}^{(n)}E_{21} & 0 & I+S_{21}^{(n)}E_{12}+S_{22}^{(n)}E_{22}
\end{array}\right)\]
and hence
$\det(I+T_n^{-1}(a)E_n)$ is equal to
\begin{equation}\det
\left(\begin{array}{cc}
I+S_{11}^{(n)}E_{11}+S_{12}^{(n)}E_{21} & S_{11}^{(n)}E_{12}+S_{12}^{(n)}E_{22}\\
S_{21}^{(n)}E_{11}+S_{22}^{(n)}E_{21} & I+S_{21}^{(n)}E_{12}+S_{22}^{(n)}E_{22}
\end{array}\right).\label{2.2}
\end{equation}
This goes to the limit
\[\det\left(\begin{array}{cc}
I+S_{11}E_{11}& S_{11}E_{12}\\
\ti{S}_{11}^\top E_{21} & I+\ti{S}_{11}^\top E_{22}
\end{array}\right)
= \det\left[\left(\begin{array}{cc}
I & 0\\ 0 & I\end{array}\right)+\left(\begin{array}{cc}
S_{11} & 0\\
0 & \ti{S}_{11}^\top
\end{array}\right)
\left(\begin{array}{cc}
E_{11} & E_{12}\\
E_{21} & E_{22}\end{array}\right)\right].\]
The assertion is now straightforward from (\ref{2.1}). $\;\:\square$

\vsg
The curve $a^\sharp$ has a natural orientation.
Under the assumption of Theorem \ref{Theo 2.1}, we may associate an argument to each point
of $a^\sharp$ such that this argument changes continuously as the point moves along the curve.
The restriction of this argument to the points in the range of $a$ defines an argument and
thus a logarithm $\log a$ of $a$. Note that if $a$ itself is continuous, then $\log a$ is
also a continuous function on the unit circle. Let $(\log a)_k$ denote the $k$th
Fourier coefficient of $\log a$. The geometric mean of $a$ is defined by
\begin{equation}
G(a)=\exp (\log a)_0=\exp \left(\frac{1}{2\pi}\int_{-\pi}^\pi \log a(e^{i\tht})\,d\tht\right).\label{GGa}
\end{equation}
It is well known that the $(1,1)$ entry of $T^{-1}(a)$ is just $1/G(a)$; see, e.g., \cite[Prop. 10.6(b)]{BoSiBu}.

\begin{exa} \label{Exa 2.2}
{\rm Suppose $m_0=1$, that is, suppose $T_n(a)$ has at most perturbations by four scalars $E_{jk}$ in its four corners.
Then
$S_{11}=\ti{S}_{11}^\top=c_{11}= 1/G(a)$
and the right-hand side of (\ref{limit}) becomes
\begin{equation}
\det\left[\left(\begin{array}{cc}
1 & 0\\ 0 & 1\end{array}\right)+\frac{1}{G(a)}
\left(\begin{array}{cc}
E_{11} & E_{12}\\
E_{21} & E_{22}\end{array}\right)\right].\label{2.3}
\end{equation}
For
\[\left(\begin{array}{cc}
E_{11} & E_{12}\\
E_{21} & E_{22}\end{array}\right)=\left(\begin{array}{cc}
1 & 0 \\
0 & 1\end{array}\right),
\quad\left(\begin{array}{cc}
0 & 1 \\
1 & 0\end{array}\right),
\quad\left(\begin{array}{cc}
1 & 1 \\
1 & 1\end{array}\right),\]
this is \[\left(1+\frac{1}{G(a)}\right)^2, \quad 1-\frac{1}{G(a)^2}, \quad \frac{2}{G(a)}+\frac{1}{G(a)^2},\] respectively.
The limit (\ref{limit}) is
zero if and only if $G(a)$ is an eigenvalue of the $2\times 2$ matrix
$-\left(\begin{array}{cc} E_{11} & E_{12} \\ E_{21} & E_{22} \end{array}\right)$. $\;\:\square$
}
\end{exa}

If the symbol $a$ is continuous, then the curve $a^\sharp$ is simply the range $a(\bT)$.
Now suppose that $a$ is sufficiently smooth, say
\begin{equation}
\sum_{k=-\iy}^\iy k^\la |a_k| <\iy,\label{sm}
\end{equation}
for some $\la >0$. The set of all $a$ satisfying (\ref{sm}) is a weighted Wiener algebra and will be denoted
by $W^\la$. If $\la >1/2$ and if $a$
has no zeros on the unit circle and winding number zero about the origin, then the
asymptotic behavior of the determinants $\det T_n(a)$ is described by Szeg\H{o}'s strong
limit theorem. This theorem says that
\begin{equation}
\det T_n(a)= G(a)^n E(a)(1+o(1))\label{Sz}
\end{equation}
where $G(a)$ is given by (\ref{GGa}) and $E(a)$ is defined by
\[E(a)=\exp \sum_{k=1}^\iy k (\log a)_{-k}(\log a)_k.\]
Formula (\ref{Sz}) may also be written in the form
\[\lim_{n \to \iy} \det T_n\left(\frac{a}{G(a)}\right)=E(a).\]
In other words, after appropriate normalization the determinants approach a finite and nonzero limit as their
order increases to infinity. In \cite[Corollary 10.38]{BoSiBu} it is shown that the $o(1)$ in (\ref{Sz})
is $o(1/n^{2\la-1})$.

\vsk
The following result is a refinement of Theorem \ref{Theo 2.1} for smooth symbols.

\begin{thm} \label{Theo 2.3}
Let $a \in W^\la$ with $\la > 1/2$ and suppose $a$ has no zeros on the unit circle and winding number zero
about the origin. Then
\[
\frac{\det (T_n(a)+E_n)}{\det T_n(a)}
=\det\left[\left(\begin{array}{cc}
I & 0\\ 0 & I\end{array}\right)+\left(\begin{array}{cc}
S_{11} & 0\\
0 & \ti{S}_{11}^\top
\end{array}\right)
\left(\begin{array}{cc}
E_{11} & E_{12}\\
E_{21} & E_{22}\end{array}\right)\right]+O\left(\frac{1}{n^{\la}}\right).\]
\end{thm}

{\em Proof.} We adopt the notations of the proof of Theorem \ref{Theo 2.1}. From Theorem 2.15
of \cite{BoSiUni} we see that $S_{11}^{(n)}=S_{11}+O(1/n^{\la})$ (entry-wise). It follows that
$S_{22}^{(n)}=[\ti{S}_{11}^{(n)}]^\top=S_{11}^\top+O(1/n^{\la})$. Let $\ell^2_\la$
be the weighted $\ell^2$ space of all sequences $x$ satisfying
\[\n x\n_{2,\la}:=\left(\sum_{n=1}^\iy n^{2\la}|x_n|^2\right)^{1/2} < \iy.\]
Theorem 7.25 of \cite{BoSiBu} implies that if $x \in \ell^2_\la$, then $T^{-1}(a)x \in \ell^2_\la$
and
\begin{equation}
\n T_n^{-1}(a)P_nx-T^{-1}(a)x\n_{2,\la} \to 0.\label{wei}
\end{equation}
Let $T_n^{-1}(a)=(c_{jk}^{(n)})_{j,k=1}^n$.
The $k$th column of $S_{12}^{(n)}$ is $(c_{n-m_0+1,k}^{(n)}, \ldots, c_{n,k}^{(n)})^\top$, while the last $m_0$
components of the $k$th column of $T^{-1}(a)$ are $c_{n-m_0+1,k}, \ldots, c_{n,k}$.

\vsk
Let $e_k$ be the
sequence which has $1$ in position $k$ and zeros elsewhere.
The convergence result (\ref{wei})
with $x=e_k$  shows that
\[\sum_{j=1}^{m_0} (n-m_0+j)^{2\la}|c_{n-m_0+j,k}^{(n)}-c_{n-m_0+j,k}|^2 \to 0.\]
This implies that $(n-m_0+j)^{2\la}|c_{n-m_0+j,k}^{(n)}-c_{n-m_0+j,k}|^2 \to 0$ and hence
\[c_{n-m_0+j,k}^{(n)}=c_{n-m_0+j,k}+o(1/n^\la).\]
Since $T^{-1}(a)e_k \in \ell^2_\la$, we also have $\sum_{n=1}^\iy n^{2\la}|c_{n,k}|^2 < \iy$,
which yields
\[c_{n-m_0+j,k} =o(1/n^\la).\]
Consequently, $c_{n-m_0+j,k}^{(n)}=o(1/n^{\la})$ and thus $S_{12}^{(n)}=O(1/n^\la)$. This in turn
tells us that
$S_{21}^{(n)}=[\ti{S}_{12}^{(n)}]^\top=o(1/n^{\la})$.
In summary, the determinant (\ref{2.2}) is
\[\det\left[\left(\begin{array}{cc}
I & 0\\ 0 & I\end{array}\right)+\left(\begin{array}{cc}
S_{11} & 0\\
0 & \ti{S}_{11}^\top
\end{array}\right)
\left(\begin{array}{cc}
E_{11} & E_{12}\\
E_{21} & E_{22}\end{array}\right)\right]+O\left(\frac{1}{n^{\la}}\right),\]
which completes the proof. $\;\:\square$

\begin{exa} \label{Exa 2.4}
{\rm Let $a(t)=(1-\mu t)(1-\nu/t)$ with $|\mu|<1, |\nu|<1$. The $n \times n$ versions of the matrices
\[\left(\begin{array}{cccc}
1+\mu\nu & -\nu & 0 & 0 \\
-\mu & 1+\mu\nu & -\nu & 0 \\
0 & -\mu & 1+\mu\nu & -\nu \\
0 & 0 & -\mu & 1+\mu\nu
 \end{array}\right), \quad
\left(\begin{array}{cccc}
1+\mu\nu & -\nu & 0 & 1 \\
-\mu & 1+\mu\nu & -\nu & 0 \\
0 & -\mu & 1+\mu\nu & -\nu \\
1 & 0 & -\mu & 1+\mu\nu
 \end{array}\right),
\]
are $T_n(a)$ and $T_n(a)+E_n$.
We have $G(a)=1$ and $E(a)=1/(1-\mu\nu)$, and hence
Szeg\H{o}'s strong limit theorem tells us that $\det T_n(a)$ has the limit $1/(1-\mu\nu)$.
Theorem \ref{Theo 2.3} may be applied with arbitrarily large $\la$.
Since $G(a)=1$ is an eigenvalue of
\[-\left(\begin{array}{cc} E_{11} & E_{12} \\ E_{21} & E_{22} \end{array}\right)
=\left(\begin{array}{rr} 0 & -1 \\ -1 & 0 \end{array}\right),\]
Example \ref{Exa 2.2} and Theorem \ref{Theo 2.3} predict that $\det(T_n(a)+E_n)/\det T_n(a)$
goes to zero faster than an arbitrary power of $1/n$.
In fact it is easy to compute the determinants exactly. We have
\begin{eqnarray*}
& & \det T_n(a)=\frac{1-(\mu\nu)^{n+1}}{1-\mu\nu},\\
& & \det (T_n(a)+E_n)= (1+\mu\nu)^2(\mu\nu)^{n-1}+\mu^{n-1}+\nu^{n-1}.
\end{eqnarray*}
This shows that the quotient $\det(T_n(a)+E_n)/\det T_n(a)$ actually decays exponentially fast to zero.
$\;\:\square$
}
\end{exa}

\section{The pure Fisher-Hartwig singulaity}\label{S3}

The symbol $a(t)=(1-t)^\ga(1-1/t)^\de$ is referred to as the pure Fisher-Hartwig singularity.
Here $\de$ and $\ga$ are complex numbers. We define
\begin{eqnarray*}
& & \xi_\de(t):=(1-1/t)^\de:=\sum_{k=0}^\iy (-1)^k\dbinom{\de}{k} t^{-k},\\
& & \eta_\ga(t):=(1-t)^\ga:=\sum_{k=0}^\iy (-1)^k \dbinom{\ga}{k}t^k
\end{eqnarray*}
and may then write $a=\xi_\de\eta_\ga$. Throughout what follows we assume that
the real parts of $\de$, $\ga$, and $\de+\ga$ are greater than $-1$. This guarantees that
$\xi_\de$, $\eta_\ga$, and $\xi_\de\eta_\ga$ are in $L^1$. Note that
the symbol (\ref{1.3}), which belongs to the $n \times n$ versions of matrix (\ref{T6}), is the pure
Fisher-Hartwig singularity $a=\xi_2\eta_2$.

\vsk
As shown in \cite[Lemma 6.18]{BoSiBu}, the $k$th Fourier coefficient of $\xi_\de\eta_\ga$
is
\[(-1)^k\frac{\Ga(1+\de+\ga)}{\Ga(\de+n+1)\Ga(\ga-n+1)}\]
in case neither $\de+n+1$ nor $\ga-n+1$ is a nonpositive integer and is equal to zero if
$\de+n+1$ or $\ga-n+1$ is a nonpositive integer. The determinants of $T_n(\xi_\de\eta_\ga)$
are known both exactly and asymptotically. Section 10.58 and Theorem 10.59 of \cite{BoSiBu}
tell us that
\begin{eqnarray}
\det T_n(\xi_\de\eta_\ga) & = & \frac{{\rm G}(1+\de){\rm G}(1+\ga)}{{\rm G}(1+\de+\ga)}\,
\frac{{\rm G}(n+1) {\rm G}(n+1+\de+\ga)}{{\rm G}(n+1+\de){\rm G}(n+1+\ga)}\label{Bar1}\\
& = & \frac{{\rm G}(1+\de){\rm G}(1+\ga)}{{\rm G}(1+\de+\ga)}\,n^{\de\ga}\,(1+o(1)),  \label{Bar2}
\end{eqnarray}
where ${\rm G}(z)$ is the Barnes function. We see in particular that $T_n(\xi_\de\eta_\ga)$ is
invertible for every $n \ge 1$. We write $T_n^{-1}(\xi_\de\eta_\ga)=(c_{jk}^{(n)}(\xi_\de\eta_\ga))_{j,k=1}^n$.

\begin{thm} \label{Theo 3.1} For each fixed $j$,
\begin{equation}
c_{jn}^{(n)}(\xi_\de\eta_\ga)=\frac{\Ga(j+\ga)}{\Ga(\de)\Ga(j)}n^{\de-\ga-1}\left(1+\frac{p_j(\xi_\de\eta_\ga)}{2n}
+O\left(\frac{1}{n^2}\right)\right)\label{ro}
\end{equation}
with
\[p_j(\xe)=(\de-j)(\de-j-1)+\de(\de-1)-(\de+\ga)(\de+\ga-1)-j(j-1)\]
and
\begin{equation}
c_{n-j,n}^{(n)}(\xe)=\frac{\Ga(j+\de)}{\Ga(\de)\Ga(j+1)}\left(1+\frac{q_j(\xe)}{2n}++O\left(\frac{1}{n^2}\right)\right)\label{ru}
\end{equation}
with
\[q_j(\xe)=(\ga-j)(\ga-j-1)+\de(\de-1)-(\de+\ga)(\de+\ga-1)-(j+1)j.\]
Furthermore, again for each fixed $j$,
\begin{equation}
c_{j1}^{(n)}(\xe)=c_{n-j+1,n}^{(n)}(\xi_\ga\eta_\de), \quad
c_{n-j,1}^{(n)}(\xe)=c_{j+1,n}^{(n)}(\xi_\ga\eta_\de). \label{li}
\end{equation}

\end{thm}

{\em Proof.} The key is the Duduchava-Roch formula, which can be found as Theorem 6.20 in \cite{BoSiBu};
see also equalities (7.87) and (7.88) of \cite{BoSiBu}.\footnote{The formula was obtained by Duduchava in the case $\ga+\de=0$
in his 1974 paper~\cite{Dud}. In 1984, Steffen Roch established the formula in the general case. With Roch's permission,
it was published in~\cite{BoSiJFA} for the first time. See~\cite[pp.~320--321]{BoSiBu} for more on the story.} This formula says that
\begin{equation}
T_n^{-1}(\xe)=\Ga_{\de,\ga} M_\ga T_n(\xi_{-\de})M_{\ga+\de}^{-1}T_n(\eta_{-\ga})M_\de, \label{DuRo}
\end{equation}
where $\Ga_{\de,\ga}=\Ga(1+\de)\Ga(1+\ga)/\Ga(1+\de+\ga)$, $M_\al$ stands for the diagonal matrix
\[M_\al={\rm diag}(\mu_1(\al), \ldots,\mu_{n}(\al)),\quad \mu_k(\al)=\frac{\Ga(k+\al)}{\Ga(1+\al)\Ga(k)},\]
$T_n(\xi_\de)$ is the upper-triangular Toeplitz matrix whose first row is
\[((\xi_{-\de})_0, \ldots, (\xi_{-\de})_{n-1}) \quad \mbox{with}\quad (\xi_{-\de})_k=\frac{\Ga(k+\de)}{\Ga(\de)\Ga(k+1)},\]
and $T_n(\eta_\ga)$ is the lower-triangular Toeplitz matrix with the first column
\[((\eta_{-\ga})_0, \ldots, (\eta_{-\ga})_{n-1})^\top \quad \mbox{with}\quad (\eta_{-\ga})_k=\frac{\Ga(k+\ga)}{\Ga(\ga)\Ga(k+1)}.\]
Let $e_n=(0,\ldots,0,1)^\top$. Using (\ref{DuRo}) it is easily seen that the $j$th component of the column
$T_n^{-1}(\xe)e_n$
is
\[c_{jn}^{(n)}(\xe)=\Ga_{\de,\ga} (\xi_{-\de})_{n-j}(\eta_{-\ga})_0\frac{\mu_j(\ga)\mu_n(\de)}{\mu_n(\de+\ga)}.\]
Inserting the above expressions for the pieces on the right we obtain
\begin{equation}
c_{jn}^{(n)}(\xe)=\frac{\Ga(j+\ga)}{\Ga(\de)\Ga(j)}\,\frac{\Ga(n-j+\de)\Ga(n+\de)}{\Ga(n-j+1)\Ga(n+\de+\ga)}.\label{cjn2}
\end{equation}
Stirling's formula gives
\begin{equation}
\frac{\Ga(n+\al)}{\Ga(n)}=n^\al\left(1+\frac{\al(\al-1)}{2n}+O\left(\frac{1}{n^2}\right)\right) \label{Gan}
\end{equation}
for every complex number $\al$. Fixing $j$ in (\ref{cjn2}), dividing numerator and denominator of (\ref{cjn2}) by $\Ga(n)^2$, and using
(\ref{Gan}) we arrive at (\ref{ro}). Replacing $j$ by $n-j$ in (\ref{cjn2}) we get
\[c_{n-j,n}^{(n)}=\frac{\Ga(j+\de)}{\Ga(\de)\Ga(j+1)}\,\frac{\Ga(n-j+\ga)\Ga(n+\de)}{\Ga(n-j)\Ga(n+\de+\ga)}.\]
Making again use of (\ref{Gan}), we obtain (\ref{ru}) for each fixed $j$.

\vsk
The numbers (\ref{ro}) and (\ref{ru}) are the upper and lower components of the last column of $T_n(\xe)$, that is,
of the column $x$ given by $T_n(\xe)x=e_n$. The entries in the first column of $T_n(\xe)$ are the entries of
the column $y$ defined by $T_n(\xe)y=e_1:=(1,0,\ldots,0)^\top$. With the counter-identity $W_n$ we therefore
have $W_nT_n(\xe)W_nW_ny=W_ne_1=e_n$, and since $W_nT_n(\xe)W_n=T_n(\xi_\ga\eta_\de)$, it follows that
$T_n(\xi_\ga\eta_\de)W_ny=e_n$. This proves (\ref{li}). $\;\:\square$

\begin{exa} \label{Exa liun}
{\rm The proof of Theorem \ref{Theo 2.1} shows that if the symbol $a$ is as in this theorem,
then the lower-left and upper-right entries of $T_n^{-1}(a)$ always approach zero as $n \to \iy$.
In Section \ref{S4} we will see that this also happens if $a \in L^1$, $a \ge 0$ almost everywhere on the unit circle,
and $\log a \in L^1$. However, Theorem \ref{Theo 3.1} reveals that in general the lower-left and upper-right entries of $T_n^{-1}(a)$
need not to converge to zero. Indeed, from (\ref{ro}) we infer that the upper-right entries of $T_n^{-1}(\xe)$ decay to zero
only if $\Ret \de-\Ret \ga < 1$, and combining (\ref{ro}) and (\ref{li}) we see that the lower-left entries of
$T_n^{-1}(\xe)$ go to zero only if $\Ret \ga -\Ret \de < 1$. Thus, both the lower-left and upper-right
entries converge to zero only if $|\Ret \ga-\Ret \de|<1$. Pure Fisher-Hartwig symbol are a nice tool
to get a first feeling for several phenomena concerning Toeplitz matrices and in particular
for disproving conjectures on such matrices! $\;\:\square$
}
\end{exa}

Theorem~\ref{Theo 3.1} is all we need to tackle the case $m_0=1$, that is, the case where $T_n(\xe)$ has at most four
scalar perturbations in the corners. From (\ref{2.1}) and (\ref{2.2}) we infer that if the $E_{jk}$ are scalars, then
\begin{equation}
\frac{\det (T_n(\xe)+E_n)}{\det T_n(\xe)}=\det\left[\left(\!\!\begin{array}{cc}
1 & 0\\ 0 & 1\end{array}\!\!\right)+\left(\!\!\begin{array}{cc}
\!c_{11}^{(n)}(\xe) & c_{1n}^{(n)}(\xe)\\
\!c_{n1}^{(n)}(\xe) & c_{nn}^{(n)}(\xe)
\end{array}\!\!\right)
\left(\!\!\begin{array}{cc}
E_{11} & E_{12}\\
E_{21} & E_{22}\end{array}\!\!\right)\right].\label{sca}
\end{equation}

\begin{exa} \label{Exa 3.2}
{\rm Suppose
\[\left(\begin{array}{cc}
E_{11} & E_{12}\\
E_{21} & E_{22}\end{array}\right)=\left(\begin{array}{cc}
0 & 1 \\
1 & 0\end{array}\right).\]
Then
\begin{eqnarray*}
\frac{\det (T_n(\xe)+E_n)}{\det T_n(\xe)} & = & \det\left[\left(\begin{array}{cc}
1 & 0\\ 0 & 1\end{array}\right)+\left(\begin{array}{cc}
c_{11}^{(n)}(\xe) & c_{1n}^{(n)}(\xe)\\
c_{n1}^{(n)}(\xe) & c_{nn}^{(n)}(\xe)
\end{array}\right)
\left(\begin{array}{cc}
0 & 1\\
1 & 0\end{array}\right)\right]\\
&  = & \det\left(\begin{array}{cc}
1+c_{1n}^{(n)}(\xe) & c_{11}^{(n)}(\xe) \\
c_{nn}^{(n)}(\xe)  & 1+c_{n1}^{(n)}(\xe) \end{array}\right),
\end{eqnarray*}
and by virtue of (\ref{li}), this equals
\begin{equation}
\det\left(\begin{array}{cc}
1+c_{1n}^{(n)}(\xe) & c_{nn}^{(n)}(\xi_\ga\eta_\de) \\
c_{nn}^{(n)}(\xe)  & 1+c_{1n}^{(n)}(\xi_\ga\eta_\de) \end{array}\right).\label{3.3}
\end{equation}
We take only the main term of (\ref{ro}) for $j=1$, and we take (\ref{ru}) for $j=0$,
in which case $q_0(\xe)= q_0(\xi_\ga\eta_\de)=-2\de\ga$. Then (\ref{3.3}) becomes
\begin{equation}
\det\left(\begin{array}{cc}
\!\!1+\frac{\Ga(1+\ga)}{\Ga(\de)}n^{\de-\ga-1}+O(n^{\Ret\de-\Ret\ga-2}) & 1-\frac{\de \ga}{n}+O\left(\frac{1}{n^2}\right)\\
\!\!1-\frac{\de \ga}{n}+O\left(\frac{1}{n^2}\right) & 1+\frac{\Ga(1+\de)}{\Ga(\ga)}n^{\ga-\de-1}+O(n^{\Ret\ga-\Ret\de-2})
\end{array}\right). \label{gade}
\end{equation}
This is
\[\frac{\Ga(1+\ga)}{\Ga(\de)}n^{\de-\ga-1}+O(n^{\Ret\de-\Ret\ga-2})\quad\mbox{for}\quad \Ret \de \ge \Ret\ga+1\]
and
\[\frac{\Ga(1+\ga)}{\Ga(\de)}n^{\de-\ga-1}+O\left(\frac{1}{n}\right)\quad\mbox{for}\quad \Ret \ga+1 >\Ret\de > \Ret\ga.\]
We know that $\det T_n(\xe)$ is asymptotically a constant times $n^{\de\ga}$.
It follows that $\det (T_n(\xe)+E_n)$ is asymptotically a constant times
\[n^{\de\ga}n^{\de-\ga-1}=n^{(\de-1)(\ga+1)}\]
provided $\Ret\de > \Ret\ga$. In the case where $\Ret \de < \Ret \ga$, we may pass to transposed matrices, which
does not change determinants but changes the roles of $\ga$ and $\de$ and therefore shows that
then $\det (T_n(\xe)+E_n)$ is asymptotically a constant times
\[n^{\de\ga}n^{\ga-\de-1}=n^{(\ga-1)(\de+1)}.\]
In summary, if $\de,\ga$ are positive real numbers, in which case
$\det T_n(\xe)$ grows with $n$, then

\begin{itemize}
\item $\det(T_n(\xe)+E_n)$ grows faster than $\det T_n(\xe)$ if $\de > \ga+1$ or $\de < \ga-1$,

\item $\det(T_n(\xe)+E_n)$ grows slower than $\det T_n(\xe)$ if $\ga-1 < \de < \ga+1$,

\item $\det(T_n(\xe)+E_n)$ decays to zero if $\ga < 1$ and $\de < 1$.  $\;\:\square$
\end{itemize}
}
\end{exa}

\vsk
The case $\de=\ga$ is especially nice and therefore deserves a separate treatment
by the following corollary.
We have $\xi_\al(t)\eta_\al(t)=|1-t|^{2\al}$. Recall that we require $\Ret \al >-1/2$ and
that for $\al=2$ we get the symbol (\ref{1.3}).
For a square matrix $A$, we abbreviate $\det A$ to $|A|$.

\begin{cor} \label{Cor 3.3}
If the $E_{jk}$ are scalars, then $\det (T_n(\xi_\al\eta_\al)+E_n)/\det T_n(\xi_\al\eta_\al)$ is
\[\left|\begin{array}{cc}
1+E_{11} & E_{12} \\
E_{21} & 1+E_{22}\end{array}\right|
+\frac{\al}{n}\left(E_{12}+E_{21}-\al(E_{11}+E_{22})-2 \al \left|\begin{array}{cc}
E_{11} & E_{12}\\
E_{21} & E_{22}\end{array}\right|\,\right)+O\left(\frac{1}{n^2}\right).\]
If in particular
\begin{equation}
\left(\begin{array}{cc}
E_{11} & E_{12}\\
E_{21} & E_{22}\end{array}\right)=\left(\begin{array}{cc}
0 & 1 \\
1 & 0\end{array}\right),\label{0110}
\end{equation}
then
\begin{equation}
\frac{\det (T_n(\xi_\al\eta_\al)+E_n)}{\det T_n(\xi_\al\eta_\al)}
=\frac{2\al(\al+1)}{n}+O\left(\frac{1}{n^2}\right).\label{alpha}
\end{equation}
\end{cor}

{\em Proof.} From Theorem \ref{Theo 3.1} we deduce that
\begin{equation}
c_{1n}^{(n)}(\xi_\al\eta_\al)=c_{n1}^{(n)}(\xi_\al\eta_\al)=\frac{\al}{n}+O\left(\frac{1}{n^2}\right)\label{rase1}
\end{equation}
and
\begin{equation}
c_{11}^{(n)}(\xi_\al\eta_\al)=c_{nn}^{(n)}(\xi_\al\eta_\al)=1-\frac{\al^2}{n}+O\left(\frac{1}{n^2}\right).\label{rase2}
\end{equation}
Thus, (\ref{sca}) equals
\[\left|\left(\begin{array}{cc}
1 & 0\\ 0 & 1\end{array}\right)+\left(\begin{array}{cc}
1-\frac{\al^2}{n}+O\left(\frac{1}{n^2}\right) & \frac{\al}{n}+O\left(\frac{1}{n^2}\right)\\
\frac{\al}{n}+O\left(\frac{1}{n^2}\right) & 1-\frac{\al^2}{n}+O\left(\frac{1}{n^2}\right)
\end{array}\right)
\left(\begin{array}{cc}
E_{11} & E_{12}\\
E_{21} & E_{22}\end{array}\right)\right|,\]
which can be simplified to the asserted expression. $\;\:\square$

\vsg
When restricted to the present context, Theorem 5 of \cite{RamSeg11} says that
\[c_{1n}^{(n)}(\xi_\al\eta_\al)=\frac{\al}{n}(1+o(1)), \quad c_{11}^{(n)}(\xi_\al\eta_\al)=\left(1-\frac{\al^2}{n}\right)(1+o(1)).\]
The second formula is probably misstated in \cite{RamSeg11} and should correctly read
\[c_{11}^{(n)}(\xi_\al\eta_\al)=1-\frac{\al^2}{n}(1+o(1)).\]
Clearly, these formulas are close to but nevertheless weaker than (\ref{rase1}) and (\ref{rase2}).

\begin{exa} \label{Exa 3.4}
{\rm We write $a_n \sim b_n$ if $a_n/b_n \to 1$.
Combining (\ref{Bar2}) and the corollary we see that the two corner
perturbations given by (\ref{0110}) lead to
\[\det(T_n(\xi_\al\eta_\al)+E_n)\sim\frac{{\rm G}(1+\al)^2}{{\rm G}(1+2\al)}\,2\al(\al+1)\,n^{\al^2-1}.\]
Thus, the exponent $\al^2$ is indeed lowered by $1$. If $k$ is a positive integer then ${\rm G}(k)=(k-2)!\ldots 2!1!$
with ${\rm G}(2)={\rm G}(1)=1$. We so obtain in particular
\begin{eqnarray*}
& & \det T_n(\xi_1\eta_1) \sim n, \quad \det(T_n(\xi_1\eta_1)+E_n) \sim 4,\\
& &  \det T_n(\xi_2\eta_2) \sim \frac{n^4}{12}, \quad \det(T_n(\xi_2\eta_2)+E_n) \sim n^3,\\
& &  \det T_n(\xi_3\eta_3) \sim \frac{n^9}{8640}, \quad \det(T_n(\xi_3\eta_3)+E_n) \sim \frac{n^8}{360}.
\end{eqnarray*}
We can of course also compute the determinants exactly. Formula (\ref{cjn2}) provides us with an exact
expression for $c_{jn}^{(n)}(\xi_\de\eta_\ga)$. It implies that
\[
c_{1n}^{(n)}(\xi_\al\eta_\al)=\al \frac{\Ga(n-1+\al)\Ga(n+\al)}{\Ga(n)\Ga(n+2\al)},\quad
c_{nn}^{(n)}(\xi_\al\eta_\al)=\frac{1}{\Ga(\al)} \frac{\Ga(n+\al)\Ga(n+\al)}{\Ga(n)\Ga(n+2\al)}.
\]
For $\al=1$, this gives
\[c_{1n}^{(n)}(\xi_1\eta_1)=\frac{n}{n+1}, \quad
c_{nn}^{(n)}(\xi_1\eta_1)=\frac{1}{n+1},\]
and inserting this in (\ref{3.3}) we obtain
\[\left|\begin{array}{cc}
1+\frac{1}{n+1} & \frac{n}{n+1}\\
\frac{n}{n+1} & 1+\frac{1}{n+1}\end{array}\right|=\frac{4}{n+1}.\]
Since $\det T_n(\xi_1\eta_1)=n+1$ due to (\ref{Bar1}), it follows that $\det(T_n(\xi_1\eta_1)+E_n)=4$
for all $n \ge 2$. Analogously, for $\al=2$ we have
\[c_{1n}^{(n)}(\xi_2\eta_2)=\frac{2n}{(n+2)(n+3)}, \quad
c_{nn}^{(n)}(\xi_2\eta_2)=\frac{n(n+1)}{(n+2)(n+3)}\]
and hence the determinant (\ref{3.3}) equals
\[\left|\begin{array}{cc}
1+\frac{2n}{(n+2)(n+3)} & \frac{n(n+1)}{(n+2)(n+3)}\\
\frac{n(n+1)}{(n+2)(n+3)} & 1+\frac{2n}{(n+2)(n+3)}\end{array}\right|=\frac{12(n+1)^2}{(n+2)^2(n+3)}.\]
The determinant $\det T_n(\xi_2\eta_2)$ is (\ref{adet}) by virtue of (\ref{Bar1}). Consequently,
\[\det (T_n(\xi_2\eta_2)+E_n)=\frac{(n+1)(n+2)^2(n+3)}{12}\cdot\frac{12(n+1)^2}{(n+2)^2(n+3)}=(n+1)^3\]
for $n \ge 2$. Similarly,
\[\det T_n(\xi_3\eta_3)=\frac{(n+1)(n+2)^2(n+3)^3(n+4)^2(n+5)}{8640}\]
for $ n \ge 1$ and
\[\det(T_n(\xi_3\eta_3)+E_n)=\frac{(n+1)(n+2)^2(n+3)[(n+2)^2+1][(n+2)^2+2]}{360}\]
for $n \ge 2$. $\;\:\square$
}
\end{exa}

To treat the case $m_0 \ge 2$, we need the matrices $S_{jk}^{(n)}$ in (\ref{Sjk}). Theorem \ref{Theo 3.1} provides
us with the first and last entries of the first and last columns of $T_n^{-1}(a)$. The entries in the four corners
$S_{jk}^{(n)}$ of $T_n^{-1}(a)$ can therefore be computed with the help of the Gohberg-Sementsul-Trench formula \cite{GoSe},
\cite{Trench}.
This formula says that if
\begin{equation}
\left(\begin{array}{c}
x_1 \\ \vdots \\ x_n\end{array}\right)
= \left(\begin{array}{c}
c_{11}^{(n)} \\ \vdots \\ c_{n1}^{(n)}\end{array}\right), \quad
\left(\begin{array}{c}
y_1 \\ \vdots \\ y_n\end{array}\right)
= \left(\begin{array}{c}
c_{1n}^{(n)} \\ \vdots \\ c_{nn}^{(n)}\end{array}\right)\label{columns}
\end{equation}
are the first and last columns of $T_n^{-1}(a)$ and $x_1 \neq 0$, then
\begin{eqnarray}
T_n^{-1}(a) & = & \frac{1}{x_1}\left(\begin{array}{ccc}
x_1 & & \\
\vdots & \ddots & \\
x_n & \ldots & x_1\end{array}\right)
\left(\begin{array}{ccc}
y_n & \ldots & y_1\\
& \ddots & \vdots \\
& & y_1 \end{array}\right)\nonumber \\
& & -\frac{1}{x_1}\left(\begin{array}{ccc}
y_0 & & \\
\vdots & \ddots & \\
y_{n-1} & \ldots & y_0\end{array}\right)
\left(\begin{array}{ccc}
x_{n+1} & \ldots & x_2\\
& \ddots & \vdots \\
& & x_{n+1} \end{array}\right),\label{GST}
\end{eqnarray}
where $x_{n+1}:=0$ and $y_0 :=0$. A full proof is also in \cite[p.~21]{HeiRo}.
If $S_{jk}^{(n)}$ has a limit $S_{jk}$, then (\ref{2.2}) implies
that
\begin{equation}
\lim_{n\to \iy}\frac{\det (T_n(a)+E_n)}{\det T_n(a)}
=\det\left[\left(\begin{array}{cc}
I & 0\\ 0 & I\end{array}\right)+\left(\begin{array}{cc}
S_{11} & S_{12}\\
S_{21} & S_{22}
\end{array}\right)
\left(\begin{array}{cc}
E_{11} & E_{12}\\
E_{21} & E_{22}\end{array}\right)\right].\label{limS}
\end{equation}

\begin{exa} \label{Exa 3.5}
{\rm Theorem \ref{Theo 3.1} applied to $a=\xi_\al\eta_\al$ shows that, for fixed $j$,
\begin{eqnarray}
& & c_{j1}^{(n)}(\xi_\al\eta_\al)=c_{n-j+1,n}^{(n)}(\xi_\al\eta_\al) \to c_j:=\dbinom{\al+j-2}{j-1}, \label{rase3}\\
& & c_{jn}^{(n)}(\xi_\al\eta_\al)=c_{n-j+1,1}^{(n)}(\xi_\al\eta_\al) \to 0.\label{rase4}
\end{eqnarray}
It follows that $S_{12}^{(n)}$ and $S_{21}^{(n)}$ converge to zero, and formula (\ref{GST}) implies that $S_{11}^{(n)}$
goes to
\[
S_{11}=\frac{1}{c_1}\left(\begin{array}{ccc}
c_1 & & \\
\vdots & \ddots & \\
c_{m_0} & \ldots & c_1 \end{array}\right)
\left(\begin{array}{ccc}
c_1 & \ldots & c_{m_0}\\
& \ddots & \vdots \\
& & c_1 \end{array}\right).\]
Since $T_n(\xi_\al\eta_\al)$ is symmetric, we see that $S_{22}^{(n)} \to \ti{S}_{11}$.
Thus, formula (\ref{limS}) becomes
\begin{equation}
\lim_{n\to \iy}\frac{\det (T_n(\xi_\al\eta_\al)+E_n)}{\det T_n(\xi_\al\eta_\al)}
=\det\left[\left(\begin{array}{cc}
I & 0\\ 0 & I\end{array}\right)+\left(\begin{array}{cc}
S_{11} & 0\\
0 & \ti{S}_{11}
\end{array}\right)
\left(\begin{array}{cc}
E_{11} & E_{12}\\
E_{21} & E_{22}\end{array}\right)\right].\label{rase5}
\end{equation}
If $m_0=1$, then $S_{11}=(1)$, and for the matrix (\ref{0110}) we get
\[\lim_{n\to \iy}\frac{\det (T_n(\xi_\al\eta_\al)+E_n)}{\det T_n(\xi_\al\eta_\al)}
=\det\left[\left(\begin{array}{cc}
1 & 0\\ 0 & 1\end{array}\right)+\left(\begin{array}{cc}
1 & 0\\
0 & 1
\end{array}\right)
\left(\begin{array}{cc}
0 & 1\\
1 & 0\end{array}\right)\right]=0.\]
This is correct but weaker than (\ref{alpha}). Notice that here we used only limits, whereas in order to establish (\ref{alpha}) we worked
with finer asymptotics. In the case $m_0=2$ we have
\[S_{11}=\left(\begin{array}{cc}
1 & \al \\
\al & 1+\al^2\end{array}\right), \quad
\ti{S}_{11}=\left(\begin{array}{cc}
1+\al^2 & \al \\
\al & 1\end{array}\right).\]
Theorem \ref{Theo 3.1} provides us with error terms
in (\ref{rase3}) and (\ref{rase4}) and thus with finer results in the case where the right-hand side of (\ref{rase5})
is zero. However, we will not embark on this issue here. $\;\:\square$
}
\end{exa}

\section{General Hermitian Fisher-Hartwig determinants} \label{S4}

We first embark on the general case where $a \in L^1$, $a \ge 0$ almost everywhere on $\bT$, and $\log a \in L^1$.
Fisher-Hartwig symbols are a special case and will be considered in the examples at the end of this section.
The constant $G(a)$ defined by (\ref{GGa}) is a finite and strictly positive real number. Let
\[\log a(t)=\sum_{k=-\iy}^\iy (\log a)_k t^k, \quad t \in \bT.\]
For $|z|<1$, we define
\[a_+(z)=\exp\sum_{k=1}^\iy (\log a)_k z^k\]
and
\[a_+^{-1}(z)=\exp\left(-\sum_{k=1}^\iy (\log a)_k z^k\right)=\sum_{k=0}^\iy (a_+^{-1})_k z^k.\]
Simon \cite[p. 144]{Sim} defines the Szeg\H{o} function associated with $a$ as
\[D(z)=\exp\left(\frac{1}{4\pi}\int_{-\pi}^\pi \frac{e^{i\tht}+z}{e^{i\tht}-z}\log a(e^{i\tht})\,d\tht\right)
=\exp\left(\frac{(\log a)_0}{2}+\sum_{k=1}^\iy (\log a)_k z^k\right).\]
Note that this is just the outer function whose modulus on $\bT$ is $|a|^{1/2}$.
Clearly, $a_+(z)=G(a)^{-1/2} D(z)$. Our assumptions imply that $T_n(a)$ is a positive definite (Hermitian) matrix for every $n \ge 1$.
We put $T_n^{-1}(a)=(c_{jk}^{(n)})_{j,k=1}^n$ and abbreviate $c_{j1}^{(n)}$ to $c_j^{(n)}$. Thus, $(c_1^{(n)}, \ldots, c_n^{(n)})^\top$
is the first column of $T_n^{-1}(a)$.

\begin{thm} \label{Theo 4.1}
For each fixed $j \ge 1$,
\begin{equation}
\lim_{n \to \iy} c_j^{(n)}=\frac{1}{G(a)}(a_+^{-1})_{j-1}, \quad \lim_{n\to \iy} c_{n-j+1}^{(n)}=0. \label{limc}
\end{equation}
\end{thm}

{\em Proof.}
The polynomial
\[\Phi_{n-1}(z)=\frac{1}{\ov{c}_1^{(n)}}(\ov{c}_n^{(n)}+\cdots+\ov{c}_2^{(n)}z^{n-2}+\ov{c}_1^{(n)}z^{n-1})\]
is known as the predictor polynomial of $a$. By virtue of \cite[Theorem 1.5.12]{Sim}, it is the $n-1$st monic
orthogonal polynomial on the unit circle $z=e^{i\tht}$ associated with the measure $d\mu(\tht)=\log a(e^{i\tht})\,d\tht/(2\pi)$.
Let $\n \Phi_{n-1}\n$ be its norm in $L^2(\bT,d\mu)$ and put
\[\ph_{n-1}(z)=\frac{1}{\n \Phi_{n-1}\n}\,\Phi_{n-1}(z)=\ka_{n-1}z^{n-1}+ \mbox{lower order powers}.\]
Thus, $\ph_{n-1}(z)=\ka_{n-1}\Phi_{n-1}(z)$. By \cite[Theorem 1.5.11(b)]{Sim}, we have
\[\ka_{n-1}^2=\prod_{j=0}^{n-2}\frac{1}{1-|\al_j|^2}=\frac{\det T_{n-1}(a)}{\det T_n(a)}=c_1^{(n)},\]
where $\al_0, \al_1, \ldots$
are the Verblunsky coefficients, and Szeg\H{o}'s theorem \cite[Theorem 2.3.1]{Sim} says that
\[\prod_{j=0}^\iy(1-|\al_j|^2)=G(a).\] It follows that $\ka_{n-1} \to G(a)^{-1/2}$ and $c_1^{(n)} \to 1/G(a)$.
By \cite[Theorem 2.4.1(iv)]{Sim}, the polynomials
\[\ph_{n-1}^*(z)=z^{n-1}\ov{\ph_{n-1}(1/\ov{z})}=\frac{\ka_{n-1}}{c_1^{(n)}}(c_1^{(n)}+\cdots+c_n^{(n)}z^{n-1})\]
converge uniformly on compact subsets of the unit disk $|z|<1$ to the function $D(z)^{-1}=G(a)^{-1/2}a_+^{-1}(z)$.
This implies that the coefficient of $z^{j-1}$ in $\ph_{n-1}^*(z)$ converges to the coefficient of $z^{j-1}$ in
$D(z)^{-1}=G(a)^{-1/2}a_+^{-1}(z)$, that is,
\[\frac{\ka_{n-1}c_j^{(n)}}{c_1^{(n)}} \to \frac{1}{G(a)^{1/2}} (a_+^{-1})_{j-1}.\]
Taking into account that $\ka_{n-1} \to G(a)^{-1/2}$ and $c_1^{(n)} \to 1/G(a)$, we conclude that $c_j^{(n)} \to
(a_+^{-1})_{j-1}/G(a)$.

\vsk
To prove the second equality of (\ref{limc}), we employ the Szeg\H{o} recursion
\[\Phi_n(z)=z\Phi_{n-1}(z)-\ov{\al}_{n-1} \Phi_{n-1}^*(z);\]
see \cite[Theorem 1.5.2]{Sim}. Written out this reads
\begin{eqnarray*}
& & \frac{1}{\ov{c}_1^{(n+1)}}(\ov{c}_{n+1}^{(n+1)}+\cdots+\ov{c}_1^{(n+1)}z^n)\\
& & = \frac{z}{\ov{c}_1^{(n)}}(\ov{c}_{n}^{(n)}+\cdots+\ov{c}_1^{(n)}z^{n-1})
- \frac{\ov{\al}_{n-1}}{{c}_1^{(n)}}({c}_{1}^{(n)}+\cdots+{c}_n^{(n)}z^{n-1})
\end{eqnarray*}
Comparing the coefficients of $z^0$ we obtain
\[\frac{\ov{c}_{n+1}^{(n+1)}}{\ov{c}_1^{(n+1)}}=-\ov{\al}_{n-1},\]
and since ${c}_1^{(n+1)} \to 1/{G(a)}$ and $\al_{n-1} \to 0$, we see that
${c}_{n+1}^{(n+1)} \to 0$. Comparison of the coefficients of $z$ gives
\[\frac{\ov{c}_{n}^{(n+1)}}{\ov{c}_1^{(n+1)}}=\frac{\ov{c}_{n}^{(n)}}{\ov{c}_1^{(n)}}
-\ov{\al}_{n-1}\frac{{c}_{2}^{(n)}}{{c}_1^{(n)}},\]
and as $c_1^{(n)} \to 1/G(a)$, $c_2^{(n)}\to (a_+^{-1})_1/G(a)$, $\al_{n-1} \to 0$,
and, by what was just proved, $c_n^{(n)} \to 0$, we arrive at the conclusion that
$c_n^{(n+1)} \to 0$. Proceeding in this way we  successively see that $c_{n-1}^{(n+1)} \to 0$,
$c_{n-2}^{(n+1)} \to 0$, etc. This proves the second assertion in (\ref{limc}). $\;\:\square$

\begin{cor} \label{Cor 4.2}
Put
\[
S_{11}=\frac{1}{c_1}\left(\begin{array}{ccc}
c_1 & & \\
\vdots & \ddots & \\
c_{m_0} & \ldots & c_1 \end{array}\right)
\left(\begin{array}{ccc}
\ov{c}_1 & \ldots & \ov{c}_{m_0}\\
& \ddots & \vdots \\
& & \ov{c}_1 \end{array}\right)\quad \mbox{with}\quad c_j=\frac{1}{G(a)}(a_+^{-1})_{j-1}.\]
Then
\[\lim_{n\to \iy}\frac{\det (T_n(a)+E_n)}{\det T_n(a)}
=\det\left[\left(\begin{array}{cc}
I & 0\\ 0 & I\end{array}\right)+\left(\begin{array}{cc}
S_{11} & 0\\
0 & \ti{S}_{11}^\top
\end{array}\right)
\left(\begin{array}{cc}
E_{11} & E_{12}\\
E_{21} & E_{22}\end{array}\right)\right].\]
\end{cor}

{\em Proof.}
Since $T_n(a)$ is Hermitian, the columns (\ref{columns}) are
\[\left(\begin{array}{c}
x_1 \\ \vdots \\ x_n\end{array}\right)
= \left(\begin{array}{c}
c_{1}^{(n)} \\ \vdots \\ c_{n}^{(n)}\end{array}\right), \quad
\left(\begin{array}{c}
y_1 \\ \vdots \\ y_n\end{array}\right)
= \left(\begin{array}{c}
\ov{c}_{n}^{(n)} \\ \vdots \\ \ov{c}_{1}^{(n)}\end{array}\right).
\]
Combining Theorem \ref{Theo 4.1} and formula (\ref{GST}) we see that
\begin{equation}
[S_{11}^{(n)} \to S_{11}, \quad S_{12}^{(n)} \to 0, \quad S_{21}^{(n)} \to 0, \quad S_{22}^{(n)}=[\ti{S}_{11}^{(n)}]^\top \to \ti{S}_{11}^\top.\label{SsS}
\end{equation}
The assertion is therefore immediate from (\ref{limS}). $\;\:\square$

\vsg
In \cite[p. 690]{DiBen} and \cite[Lemma 3.2]{NouVas} it is shown that if $a$ is a (real-valued and nonnegative)
trigonometric polynomial, then the  norms of $S_{11}^{(n)}, S_{12}^{(n)},  S_{21}^{(n)}, S_{22}^{(n)}$ remain bounded
as $n \to \iy$. From~(\ref{SsS}) we see that, under the sole assumption
that $a \in L^1$, $a \ge 0$ almost everywhere on $\bT$, and $\log a \in L^1$,
these matrices even converge to limits.

\vsk
The following two examples concern perturbations of Hermitian Fisher-Hartwig matrices.

\begin{exa} \label{Exa 4.3}
{\rm Let $a(t)=\xi_\al(t)\eta_{\al}(t)b(t)=|1-t|^{2\al}b(t)$ where $\al > -1/2$ is a real number and
$b$ is a twice continuously differentiable and strictly positive function on the unit circle. Then
\[\det T_n(a) \sim G(b)^n n^{\al^2} E_*(a)\]
with some nonzero constant $E_*(a)$; see \cite[Lemma 6.47]{BoSiAkad} and \cite[Theorem 5.44]{BoSiUni}.
In this case Corollary~\ref{Cor 4.2} is applicable. We have
$c_j=(\eta_{-\al}b_+^{-1})_{j-1}/G(b)$ and hence
\begin{eqnarray*}
& & c_1=1,\\
& & c_2=(b_+^{-1})_1+\al,\\
& & c_3= (b_+^{-1})_2+(b_+^{-1})_1 \al+ \al(\al+1)/2,\\
& & \mbox{and so forth.}
\end{eqnarray*}
For the pure singularity, i.e., when $b(t)$ is identically $1$, we get
\[c_1=1, \quad c_2=\al, \quad c_3=\al(\al+1)/2,\]
and $S_{11}$ takes the same form as in Example~\ref{Exa 3.5}. $\;\:\square$
}
\end{exa}

\begin{exa} \label{Exa 4.4}
{\rm Now suppose
\[a(t)=|t_1-t|^{2\al_1} \cdots |t_r-t|^{2\al_r} b(t)\]
where $t_j$ are distinct points on $\bT$, $\al_j$ are real numbers in $(-1/2,1/2)$, and $b$
is a twice continuously differentiable and strictly positive function on $\bT$. This time
\[\det T_n(a)=G(b)^n n^{\al_1^2+\cdots+\al_r^2} E_{**}(a)\]
with some nonzero constant $E_{**}(a)$; see \cite[Theorem 5.47]{BoSiUni}. Corollary~\ref{Cor 4.2}
is again applicable. If, for example, $a(t)=|t_1-t|^{2\al_1}|t_2-t|^{2\al_2}$, then
\begin{eqnarray*}
& & c_1=1,\\
& & c_2=\frac{\al_1}{t_1}+\frac{\al_2}{t_2},\\
& & c_3=\frac{\al_1(\al_1+1)}{2t_1^2}+\frac{\al_1\al_2}{t_1t_2}+\frac{\al_2(\al_2+1)}{2t_2^2}. \quad \square
\end{eqnarray*}
}
\end{exa}

The values for $c_j$ given in Example \ref{Exa 4.3} can also be derived from \cite[Lemma 1]{RamSeg11}.
Moreover, Theorem 5 of \cite{RamSeg11}, with the surmised correction mentioned above after Corollary~\ref{Cor 3.3},
gives the second term in the asymptotics of $c_j^{(n)}$ for symbols as in Example~\ref{Exa 4.3}.
In the case of two singularities with the same exponent, that is, for $a(t)=|t_1-t|^{2\al} |t_2-t|^{2\al} b(t)$
with $-1/2 < \al < 1/2$,  which is a special case of Example~\ref{Exa 4.4},
Theorem 7 of \cite{Ram13} says that $c_j^{(n)}=(a_+^{-1})_{j-1}/G(a)+O(1/n)$, which
is stronger than our result $c_j^{(n)}=(a_+^{-1})_{j-1}/G(a)+o(1)$.

\bigskip
A. B\"ottcher, Fakult\"at f\"ur Mathematik, TU Chemnitz, 09107 Chemnitz, Germany

{\tt aboettch@mathematik.tu-chemnitz.de}

\medskip
L. Fukshansky, Department of Mathematics,  Claremont McKenna College,

850 Columbia Ave,
Claremont, CA 91711, USA

{\tt lenny@cmc.edu}

\medskip
S. R. Garcia, Department of Mathematics, Pomona College,

610 N. College Ave, Claremont, CA 91711, USA

{\tt stephan.garcia@pomona.edu}, URL: \url{http://pages.pomona.edu/~sg064747/}

\medskip
H. Maharaj, 8543 Hillside Road, Rancho Cucamonga, CA 91701, USA

{\tt hmahara@g.clemson.edu}

\end{document}